\documentclass{amsart}
\usepackage{amssymb}
\usepackage{amsfonts}
\usepackage{amsmath}

\newtheorem{theorem}{Theorem}
\theoremstyle{plain}

\newtheorem{corollary}{Corollary}
\newtheorem{assumption}{Assumption}

\newtheorem{lemma}{Lemma}

\newtheorem{proposition}{Proposition}

\numberwithin{equation}{section}
\newcommand{\Prob}{\mathbb{P}} 
\newcommand{\B}{\mathfrak{B}} 
\newcommand{\R}{\mathbb{R}} 
\newcommand{\N}{\mathbb{N}} 
\newcommand{\RR}{\overline{\mathbb{R}}}
\newcommand{\restr}[2]{#1\vert #2} 
\newcommand{\lrestr}[2]{\left.#1\right\vert #2}

\newcommand{\abs}[1]{\vert #1\vert} 
\newcommand{\dabs}[1]{\left\vert #1\right\vert} 
\newcommand{\norm}[1]{\Vert #1\Vert} 
 
\newcommand{\cond}[3]{#1(#2\vert #3)} 
\newcommand{\rcond}[3]{#1\left(#2\left\vert #3\right.\right)} 
 
\newcommand{\condP}[2]{P(#1\vert #2)}

\newcommand{\neta}[1]{\net{#1}{\alpha}{A}} 
\newcommand{\nneta}[1]{\nnet{#1}{\alpha}{A}} 
\newcommand{\net}[3]{\left\langle #1_{#2}\right\rangle_{#2\in #3} } 
\newcommand{\nnet}[3]{\left\langle #1\right\rangle_{#2\in #3} } 
\newcommand{\seq}[2]{\net{#1}{#2}{\mathbb{N}}} 
\newcommand{\sseq}[2]{\nnet{#1}{#2}{\mathbb{N}}} 
\newcommand{\seqn}[1]{\seq{#1}{n}} 
\newcommand{\sseqn}[1]{\sseq{#1}{n}}

\newcommand{\A}{\mathcal{A}} 
\newcommand{\G}{\mathcal{G}} 
\newcommand{\F}{\mathcal{F}} 
\newcommand{\U}{\mathcal{U}} 
 
\newcommand{\C}{\mathcal{C}} 
\newcommand{\K}{\mathcal{K}}
\newcommand{\Pd}{\mathcal{P}^d}  
\newcommand{\Pn}{\mathcal{P}^{d_n}}  
\newcommand{\AP}{\mathfrak{A}(P)} 
\newcommand{\OR}{\bar{\Omega}} 
\newcommand{\set} [1]{\mathbf{1}_{#1}} 
\newcommand{\Su}{\mathfrak{S}} 
\newcommand{\M}{\mathfrak{M}} 
\newcommand{\I}[1]{\mathcal{I}_{#1}} 
 
\newcommand{\Icond}[3]{#1(\left.#2\right\vert\I{#3})}  
\newcommand{\Idcond}[3]{#1\left(#2\left\vert\I{#3}\right.\right)}

\begin{document}

\title{Finitely Additive Supermartingales}
\author{Gianluca Cassese} 
\address{Universit\`{a} del Salento and University of Lugano}
\email{g.cassese@economia.unile.it}
\curraddr{Dipartimento di Scienze Economiche e Matematico-Statistiche,
Ecotekne, via per Monteroni, 73100 Lecce}
\thanks{I am in debt with an anonymous referee for several helping suggestions. 
All remaining erros are my own.}
\date
\today


\maketitle

\begin{abstract}
The concept of finitely additive supermartingales, originally due to
Bochner, is revived and developed. We exploit it to study measure
decompositions over filtered probability spaces and the properties of the
associated Dol\'{e}ans-Dade measure. We obtain versions of the Doob Meyer
decomposition and, as an application, we establish a version of the
Bichteler and Dellacherie theorem with no exogenous probability measure.
\end{abstract}

\section{Introduction}

In the classical theory of probability one often encounters situations in which
countable additivity fails. Broadly speaking, these fall into two main classes 
of problems: those involving duality on the space $L^\infty$ and those in which 
the underlying $\sigma$ algebra needs to be extended, e.g. to overcome the lack 
of measurability of some random quantity. Both situtations are well documented 
in applications. \cite{karatzas} is a recent example of the former kind of situa%
tions arising in mathematical finance. In the area of weak convergence it is
well known that even the classical empirical process is not Borel measurable
in the space $D[0,1]$ when the latter is equipped with the non separable 
topology induced by the supremum norm. Dudley \cite{dudley} illustrates a 
number of situations relevant for empirical processes in which measurability fails.
To overcome these drawbacks, a new approach, based on outer expectation, was 
developed by Hoffmann-Jorgensen and, in a systematic way, in the book by van 
der Vaart and Wellner \cite{wellner}. More recently Berti and Rigo \cite{berti rigo} 
have shown that such notion of weak convergence has an exact translation in the 
language of finite additivity. 

In this paper we fix an algebra $\A$ of subsets of some set $\Omega$ and an 
increasing family $(\A_t:t\in\R_+)$ of sub algebras of $\A$, a filtration ($\F$ 
and $\F_t$ will hereafter denote the $\sigma$ algebras generated by $\A$ and 
$\A_t$ respectively). To illustrate our topic, consider the quantity 
$m(f_t)$ where $f=(f_t:t\in\R_+)$ is an adapted process and $m$ a finitely 
additive probability on $\A$, that is $m$ is a positive, finitely additive set 
function on $\A$ (in symbols $m\in ba(\A)_+$) and $m(\Omega)=1$. When dealing with 
finitely additive expectation, of special importance are the structural properties of 
$m$ such as decompositions, particularly the one of Yosida and Hewitt \cite{yosida 
and hewitt}. In our setting, however, what matters are the properties of $m$ 
\textit{conditional} on $\A_t$ and the focus then shifts from the finitely additive 
measure $m$ to the finitely additive process $(m_t:t\in\R_+)$ where $m_t=\restr{m}{\A_t}$; 
or even $(m_t^c:t\in\R_+)$, where $m_t^c$ and $m_t^\perp$ designate the countably and 
purely finitely additive components of $m_t$ (in the sequel the spaces of countably 
and purely finitely additive set functions on an algebra $\G$ will be indicated with 
the symbols $ca(\G)$ and $pfa(\G)$ respectively while $\Prob(\G)$ will be used to 
designate full probabilities, i.e. countably additive, on $\G$). The inclusion 
$\A_t\subset\A_u$ for $u\geq t$ implies $\restr{m_u^c}{\A_t}\leq m_t$ and 
$m_t^\perp\perp\restr{m_u^c}{A_t}$ i.e. $\restr{m_u^c}{\A_t}=m_t\wedge(\restr{m_u^c}{\A_t})%
\leq m_t^c+\left(m_t^\perp\wedge(\restr{m_u^c}{\A_t})\right)=m_t^c$, a conclusion which extends to 
any decomposition $m_t=m_t^e+m_t^p$, such that $m_t^p\perp\restr{m_u^e}{\A_t}$ -- see e.g. 
Lemma \ref{lemma measure decomposition}. $(m_t:t\in\R_+)$ and $(m^c_t:t\in\R_+)$ turn 
thus out being finitely additive supermartingales, a concept introduced by Bochner, 
in a number of little known papers -- \cite{bochner 47}, \cite{bochner 50} and \cite{bochner} 
-- and later revived by Armstrong in \cite{armstrong 83} and \cite{armstrong 85}.

More formally, a finitely additive stochastic process $\xi=(\xi_t:t\in\R_+)$ 
is an element of the vector lattice $\prod_{t\in\R_+}ba(\A_t)$ endowed with the 
order induced by each coordinate space. $\xi$ is a finitely additive supermartingale if
\begin{equation}
\label{supermartingale}
\xi _t(F)\geq\xi_u(F)\qquad\qquad F\in\A_t,t\leq u 
\end{equation}
The symbol $\Su$ designates the set of finitely additive supermartingales such that
$\norm{\xi}\equiv\sup_{t\in\R_+}\norm{\xi _t}<\infty$.
With no loss of generality we put $\A_0=\bigcap_{t\in\R_+}\A_t$ and 
$\A=\bigcup_{t\in\R_+}\A_t$ and define, for $F\in\A$ and $G\in\A_0$, 
$\xi_\infty(F)=\inf_{\{t:F\in\A_t\}}\xi_t(F)$ and $\xi_0(G)=\sup_t\xi_t(G)$
a choice that will allow us to replace $\R_+$ by $\RR_+=\R_+\cup\{\infty\}$ 
when necessary. We also use the symbols $\OR=\Omega\times\RR_+$ and 
$\bar\F\equiv\F\otimes\{\varnothing,\R_+\}$.

We use repeatedly the following corollary of the Hahn-Banach theorem 
(see e.g. \cite[3.2.3(b) and 3.2.10]{rao}).
\begin{lemma}
\label{lemma rao}
If $\Sigma_0\subset\Sigma$ are algebras of subsets of some set $S$ and 
$\mu_0\in ba(\Sigma_0)_+$ then there exists $\mu\in ba(\Sigma)_+$ such that 
$\restr{\mu}{\Sigma_0}=\mu_0$.
\end{lemma} 

Although all processes in this paper are indexed by $\R_+$, we often do not use
but the order properties of the real numbers so that some of the results
that follow carry over almost unchanged to the case of a linearly ordered index set.

\section{Finitely Additive Conditional Expectation}
\label{sec conditioning}

The absence of a satisfactory concept of conditional expectation in the
finitely additive setting, a major argument in favour of countable 
additivity, is a direct consequence of the failure of the Radon Nikodym theorem. 
The operator defined hereafter, e.g., provides an extension of such fundamental 
concept which is suitable for many analytical purposes but lacks some of the 
properties which matter for the sake of its statistical interpretation (a different 
proof of the following result appears in \cite[proposition 2.1, p. 27]{io MAFI}).

\begin{theorem}
\label{theorem conditioning}
Let $\mathcal{H}$ be an algebra of subsets of $\Omega$, $\G$ a sub $\sigma$ 
algebra of $\mathcal{H}$, and $\mu\in ba(\mathcal{H})_+$. 
Let $\restr{\mu}{\G}$ decompose as $\lambda+\pi$, with $\lambda\in ca(\G)_+$, 
$\pi\in ba(\G)_+$ and $\lambda\perp\pi$ and define
\begin{equation}
\I{\pi}=\{G\in\G:\pi(G)=0\}
\label{I}
\end{equation}%
If $f\in L^1(\mu)$ there exists a unique $\Icond{\mu}{f}{\pi}\in L^1(\lambda)$
such that 
\begin{equation}
\mu(f\set{I})=\mu(\Icond{\mu}{f}{\pi}\set{I})=\lambda(\Icond{\mu}{f}{\pi}\set{I})\qquad I\in\I{\pi }
\label{conditioning}
\end{equation}
and
\begin{equation}
\Icond{\mu}{f\set{G}}{\pi}=\Icond{\mu}{f}{\pi} \set{G}\qquad G\in\G
\label{take out}
\end{equation}%
$\Icond{\mu}{\cdot}{\pi}:L^1(\mu)\rightarrow L^1(\lambda)$ is a positive, 
linear operator with $\norm{\Icond{\mu}{\cdot}{\pi }} =1$.
\end{theorem}

\begin{proof}
Being closed with respect to finite unions, $\I{\pi }$ is a
directed set relatively to inclusion. Since $\lambda\perp\pi$ and $%
\G$ is a $\sigma$ algebra, for each $\epsilon >0$ there exists $%
I\in\I{\pi }$ such that $\lambda(I^c)\leq\epsilon$: i.e. 
$\lambda(G)=\lim_{I\in\I{\pi}}\lambda(IG)$, $G\in\G$. Let 
$f\in L^1(\mu) $. Any solution $p(f)\in L^1(\lambda)$ to (\ref{conditioning}) 
must then satisfy 
$$
\lambda(p(f)\set{G})=\lim_{I\in\mathcal{I}_{\pi}}\lambda(p(f)\set{G\cap I}) 
=\lim_{I\in\mathcal{I}_\pi}\mu(f\set{G\cap I})\qquad G\in\G
$$
and is therefore unique $P$ a.s.: by considering $f^+$ and $f^-$
separately we can (and will) thus restrict to the case in which $f\in L^1(\mu)_+$.

Let $\mu_f\in ba(\G)_+$ be defined implicitly by letting 
\begin{equation}
\mu_f(G)=\lim_{I\in\I{\pi}}\mu(f\set{G\cap I})\qquad G\in\G 
\label{uf}
\end{equation}%
The limit in (\ref{uf}) exists uniformly with respect to $G\in\G$.
In fact, $I\in\I{\pi}$ implies $\mu_f(I)=\mu(f\set{I})$ and 
$\lim_{I\in\I{\pi}}\mu_f(I^c) =0$ so that 
$$
0\leq\mu _f(G) -\mu(f\set{G\cap I}) =\mu_f(G)-\mu_f(G\cap I)\leq\mu_f(I^c)
$$
Let $\seqn{G}\subset\G$ be such that $\lim_n\lambda(G_n) =0$ and $I\in\I{\pi}$. 
Then $\lim_n\mu(G_n\cap I)=\lim_n\lambda(G_n\cap I)=0$ so that $\lim_n\mu(f\set{G_n\cap I})=0$, 
by absolute continuity of the finitely additive integral \cite[III.2.20(b)]{bible}. 
Then \cite[I.7.6]{bible} 
$$
\lim_n\mu_f(G_n)=\lim_n\lim_{I\in\mathcal{I}_{\pi}}\mu(f\set{G_n\cap I})=
\lim_{I\in\I{\pi}}\lim_n\mu(f\set{G_n\cap I}) =0
$$
i.e. $\mu _f\ll\lambda$. (\ref{conditioning}) follows by letting $\Icond{\mu}{f}{\pi}
\in L^1(\lambda)_+$ be the corresponding Radon Nikodym derivative; (\ref{take out}) 
from $IG\in\I{\pi}$ whenever $I\in\I{\pi}$ and $G\in\G$. 
$\Icond{\mu}{\cdot}{\pi}$ is linear and positive as $\mu$ is. If $f\in L^1(\mu)$,
$$
\lambda(\abs{\Icond{\mu}{f}{\pi}})\leq
\lim_{I\in\I{\pi }}\lambda(\Icond{\mu}{\abs{f}}{\pi}\set{I})=
\lim_{I\in\I{\pi }}\mu(\abs{f} \set{I})\leq\norm{f}
$$
with equality if $f$ is the indicator of some $I\in\I{\pi }$ i.e. 
$\norm{\Icond{\mu}{\cdot}{\pi}}=1$.
\end{proof}

Referring to $\Icond{\mu}{\cdot}{\pi}$ as \textquotedblleft conditional 
expectation\textquotedblright\ is just a convenient abuse of terminology 
as the law of total probability $\mu(f) =\mu(\Icond{\mu}{f}{\pi})$, which 
is at the basis of the statistical interpretation of this concept \cite%
[p. 1229]{kadane}, will in general not hold unless $\restr{\mu}{\G}\in ca(\G)$%
\footnote{The failure of this property for reasonable definitions of finitely 
additive conditional expectation is well known since the work of Dubins \cite%
{dubins}}. Of course, if $\mu\in ca(\mathcal{H})$ the above concept of conditional 
expectation would coincide (by uniqueness) with the traditional one.

\section{The Dol\'{e}ans-Dade Measure}
\label{sec Doleans}

In the early works of Dol\'{e}ans-Dade \cite{doleans-dade}, F\"{o}llmer 
\cite{follmer} and Metivier and Pellaumail \cite{metivier}, supermartingales 
were associated with measures over predictable rectangles. We address this 
issue in the present setting. The claims and the proofs of this this section 
remain true if we replace $\R_+$ by any linearly ordered index set.

Denote by $\mathcal{R}$ the collection of all sets of the form 
\begin{equation}
\label{R}
F_0\times\{0\}\cup\bigcup_{n=1}^NF_n\times]t_n,\infty[  
\end{equation}%
where $F_0\in\A_0$, $N\in\N$ and, for each $N\geq n>m\geq 1$, 
$F_n\in\A_{t_n}$ and $F_n\cap F_m=\varnothing $. $\mathcal{R}$ is closed 
with respect to intersection and contains $\OR$ and 
$\varnothing $. We denote by $\mathcal{P}$ the algebra generated by 
$\mathcal{R}$: each $F\in\mathcal{P}$ takes then the form of a disjoint union
\begin{equation}
\label{P}
F_0\times\{0\}\cup\bigcup_{n=1}^NF_n\times]t_n,u_n]  
\end{equation}
with $F_0\in\A_0$, $t_n,u_n\in\R_+ $, $F_n\in\A_{t_n}$. We also denote by 
$\bar{\mathcal{P}}$ the collection defined as in (\ref{P}) but with $\A_t$
replaced by $\A$ for each $t\in\R_+$. Let 
\begin{equation}
\label{M}
\M=\left\{\bar x\in ba\left(2^{\OR}\right)_+:\lim_t\bar x\left(\Omega\times]t,\infty[\right)=0\right\}
\end{equation}

\begin{theorem}
\label{theorem extension}
$\xi\in\Su$ if an only if there exists $\bar x\in \M$ and $\lambda\in ba(\A_\infty)$
such that
\begin{equation}
\label{Doleans}
\xi_t(F)=\lambda(F)+\bar x(F\times]t,\infty[)\qquad t\in\R_+,\ F\in\A_t
\end{equation}
\end{theorem}
\begin{proof}
Assume that $\xi\in\Su$ and, replacing $\xi_t$ with $\xi_t-\xi_\infty$, assume
also that $\xi_\infty=0$. For each
$F_0\times\{0\}\cup\bigcup_{n=1}^NF_n\times]t_n,\infty[\in\mathcal{R}$ define
the quantity
\begin{equation}
x\left(F_0\times\{0\}\cup\bigcup_{n=1}^NF_n\times]t_n,\infty[\right)=
\sum_{n=1}^N\xi_{t_n}(F_n)  
\label{x}
\end{equation}%
Let $F=F_0\times\{0\}\cup\bigcup_{n=1}^NF_n\times]t_n,\infty[$ and 
$G=G_0\times\{0\}\cup\bigcup_{k=1}^KG_k\times]u_k,\infty[$  
be sets in $\mathcal{R}$. To prove that $\mathcal{R}$ is a lattice, write 
$F_0'=F_0\cap G_0^c$, $G_0'=G_0$ and, for $n,k>0$, 
$F_n'=F_n\cap\bigcap_{\{k>0:u_k\leq t_n\}}G_k^c$ 
and $G_k'=G_k\cap\bigcap_{\{n>0:t_n<u_k\} }F_n^c$. We obtain 
$$
F\cup G=
(F_0'\cup G_0')\times\{0\} 
\cup\left(\bigcup_{n=1}^NF_n'\times]t_n,\infty[\right)
\cup\left(\bigcup_{k=1}^KG_k'\times]u_k,\infty[\right)
\in\mathcal{R}
$$
Rearrange the collection $\{t_n,u_k:1\leq n\leq N,1\leq\,k\leq K\}$ as 
$\left\langle\gamma_i\right\rangle_{i=1}^I$ with $\gamma_i\geq\gamma _{i+1}$ 
and set conventionally $\xi_{\gamma_0}=0$. For $1\leq i\leq I$, let $\hat{\psi}_i$ be 
a positive extension of 
$\restr{(\xi_{\gamma_i}-\xi _{\gamma _{i-1}})}{\A_{\gamma_i}}$ 
to $\A$ and $\hat{\xi}_{\gamma_i}=\sum_{j=1}^i\hat{\psi}_j$. 
Then, $\restr{\hat{\xi}_{\gamma_i}}{\A_{\gamma_i}}=\xi_{\gamma_i}$ and 
$\hat{\xi}_{\gamma_{i+1}}\geq\hat{\xi}_{\gamma_i}\geq 0$, i.e. 
$\hat{\xi}_{t_n}\geq\hat{\xi}_{u_k}$ whenever $t_n\leq u_k$. 
For $0<k\leq K$, $\bigcup_{n=1}^{N}(G_k\cap F_n)\times]u_k,t_n]\subset F^c\cap G$ 
(as $F_n\cap F_m=\varnothing$ for $n>m>0$). Then $G\subset F$\ implies 
$G_k\cap F_n=\varnothing$ for all $1\leq k\leq K$ and $1\leq n\leq N$ such
that $u_k<t_n$ that is $G_k=\bigcup_{\{n:t_n\leq u_k\}}(G_k\cap F_n)$. Therefore,
\begin{equation}
\sum_{k=1}^K\xi_{u_k}(G_k)
=\sum_{\{1\leq n\leq N,1\leq k\leq K:t_n\leq u_k\}}\xi_{u_k}(G_k\cap F_n)\leq
\sum_{\{1\leq n\leq N,1\leq k\leq K:t_n\leq u_k\}}
\hat{\xi}_{t_n}(G_k\cap F_n)\leq\sum_{n=1}^N\xi_{t_n}(F_n)  
\label{major}
\end{equation}%
The set function $x$ defined in (\ref{x}) is then monotonic and \textit{a fortiori} well defined.
If $F=\varnothing$ then $\bigcup_{n=1}^NF_n=\varnothing$ so that $x(F) =0$; moreover 
\begin{eqnarray*}
x(F\cup G) &=&\sum_{n=1}^N\xi _{t_n}(F_n')+\sum_{k=1}^K\xi_{u_k}(G_k')\\
 &=&\sum_{n=1}^N\xi_{t_n}(F_n)-\sum_{\{1\leq n\leq N,1\leq k\leq K:t_n\geq u_k\}}\xi_{t_n}(F_n\cap G_k) \\
&&+\sum_{k=1}^K\xi_{u_k}(G_k)-\sum_{\{1\leq n\leq N,1\leq k\leq K:t_n<u_k\} }\xi_{u_k}(F_n\cap G_k)\\
 &=&\sum_{n=1}^N\xi_{t_n}(F_n)+\sum_{k=1}^K\xi_{u_k}(G_k)-\sum_{k,n\geq 1}\xi_{t_n\vee u_k}(F_n\cap G_k)\\
 &=&x(F)+x(G)-x(F\cap G)
\end{eqnarray*}
In other words, $x$ is a strongly additive set function on a lattice of sets
which contains $\varnothing$ as well as $\OR$: as such 
\cite[3.1.6, 3.2.1 and 3.2.5]{rao} it admits an extension $\hat{x}$ to $%
\mathcal{P}$. If $F\in\mathcal{P}$ is as in (\ref{P}), then%
$$
\hat{x}(F)=\sum_{n=1}^N\{x(F_n\times ]t_n,\infty[)-x(F_n\times ]u_n\vee t_n,\infty[)\}
$$
$\hat{x}$ is thus unique and, in view of (\ref{major}), positive. But then, 
$\hat{x}$ admits a positive extension $\bar x$ to $2^{\OR}$, by Lemma \ref{lemma rao}.
If $\bar x\in\M$ and $\lambda\in ba(\A_\infty)$ are given and $\xi$ is defined as in 
(\ref{Doleans}) then it is clear that $\xi\in\Su$ with $\xi_\infty=\lambda$
\end{proof}

If $\xi\in\Su$, then denote by $\M(\xi)$ the collections of those $\bar x\in\M$ 
meeting (\ref{Doleans}). Each $\bar x\in\M(\xi)$ will be referred to as
a Dol\'{e}ans-Dade measure associated to $\xi$. Disregarding the apparent 
arbitrariness implicit in the existence of a multiplicity of such measures, 
there are several remarkably simple implications of Theorem \ref{theorem extension} 
on such relevant issues as the decomposition and extension of supermartingales that 
are spelled out in the next corollaries where the following notation is used:
let $\M^{uc},\M^{up}\subset\M$ consist of measures $m$ such that 
$\restr{m}{\bar\F}\in ca(\bar\F)$ and $\restr{m}{\bar\F}\in pfa(\bar\F)$, 
respectively. We set 
\begin{equation}
\label{Su}
\Su^{uc}=\left\{\xi\in\Su:\xi_\infty\in ca(\A),\ \M(\xi)\cap\M^{uc}\neq\varnothing\right\}
\quad\text{and}\quad
\Su^{up}=\left\{\xi\in\Su:\xi_\infty\in pfa(\A),
\ \subset\M^{up}\right\}
\end{equation}
A supermartingale belonging to $\Su^{uc}$ (resp. $\Su^{up}$) will be called
uniformly countably additive (resp. uniformly purely finitely additive).

\begin{corollary}
\label{corollary ADM}
Any $\xi\in\Su$ admits a decomposition 
\begin{equation}
\label{ADM} 
\xi=\mu-\alpha
\end{equation}
where $\mu$ is a finitely additive 
martingale and $\alpha$ a positive, finitely additive increasing process (as 
defined in \cite[p. 287] {armstrong 83}). Moreover, the following are equivalent:
\begin{enumerate}
\item[(\textit{i})]$\xi\in\Su^{uc}$;
\item[(\textit{ii})] $\mu$ and $\alpha$ in (\ref{ADM}) may be chosen such that 
$\mu_\infty,\alpha_\infty\in ca(\A)$ (and thus so that $\mu$ and $\alpha$ are 
countably additive processes);
\item[(\textit{iii})] there exists $\lambda\in ca(\A)_+$ such that 
$\abs{\xi_t}\leq\restr{\lambda}{\A_t}$ for each $t\in\R_+$.
\end{enumerate}
\end{corollary}

\begin{proof}
Let $\bar{x}\in\M(\xi)$ and define
\begin{equation}
\label{mua} 
\mu_t(F)=\xi_\infty(F)+\bar x(F\times\R_+)\quad\text{and}\quad
\alpha_t(F)=\bar x(F\times [0,t]) \qquad F\in\A_t
\end{equation}
Then (\ref{ADM}) follows from (\ref{Doleans}). In fact $\mu$ is a 
finitely additive martingale while $\alpha$ extends to an increasing 
family $(\bar\alpha_t:t\in\R_+)$ of measures on $\A$ such that 
$\inf_t\norm{\bar\alpha_t}=\norm{\alpha_0}=0$. If (\textit{i}) holds, 
then upon choosing $\bar x\in\M(\xi)$ such that $\restr{\bar x}{\bar\F}$, 
(\ref{mua}) implies (\textit{ii}). If (\textit{ii}) holds let 
$\lambda=\abs{\mu_\infty}+2\alpha_\infty$. Then in restriction to $\A_t$ 
we obtain $\abs{\xi_t}\leq\abs{\xi_\infty}+(\xi_t-\xi_\infty)\leq\abs{\mu_\infty}+
\alpha_\infty+(\alpha_\infty-\alpha_t)\leq\lambda$ and (\textit{iii}) 
follows. Assume (\textit{iii}), then $\xi_\infty\in ca(\A)$. Let 
$U=\{0=t_1\leq\ldots\leq t_N\}$ and define $\bar\zeta^U_{t_N}\in ca(\A)_+$ 
to be an extension of $\xi_{t_N}-\restr{\xi_\infty}{\A_{t_N}}$ to $\A$ 
dominated by $\lambda-\xi_\infty$ and set $\bar\xi^U_{t_N}=\xi_\infty+\zeta^U_{t_N}$; 
likewise, for $n<N$ let $\bar\zeta^U_{t_n}\in ca(\A)_+$ be an extension 
of $\xi_{t_n}-\restr{\xi_{t_{n+1}}}{\A_{t_n}}$ to $\A$ dominated by 
$\lambda-\bar\xi^U_{t_{n+1}}$ and set $\bar\xi^U_{t_n}=\bar\xi^U_{t_{n+1}}+\bar\zeta^U_{t_n}$. 
Define $\bar\xi^U=\sum_{n=1}^N\bar\xi^U_{t_n}\set{[t_n,t_{n+1}[}$ a map 
from $\R_+$ to $ba(\A)$. One easily establishes that $\xi_\infty\leq\bar\xi^U_t\leq\lambda$ 
for each $t\in\R_+$, i.e. $\bar\xi^U\in[\xi_\infty,\lambda]^{\R_+}$, that 
$\bar\xi^U$ is decreasing and that $\restr{\bar\xi^U_t}{\A_t}=\xi_t$ when 
$t\in U$. If $ba(\A)^{\R_+}$ is equipped with the product topology obtained 
after endowing each coordinate space with the weak$^*$ topology, we conclude 
that $[\xi_\infty,\lambda]^{\R_+}$ is compact and that the net 
$\nnet{\bar\xi^U}{U}{\U}$, with $\U$ denoting the collection of finite subsets
of $\R_+$ directed by inclusion, admits a cluster point $\bar\xi$. Then necessarily,
$\bar\xi$ is decreasing and $\restr{\bar\xi_t}{\A_t}=\xi_t$ for each $t\in\R_+$. 
The same argument used in the proof of Theorem \ref{theorem extension} shows that 
the quantity $\sum_{n=1}^N(\bar\xi_{t_n}-\bar\xi_{u_n})(F_n)$, where 
$F_0\times\{0\}\cup\bigcup_{n=1}^NF_n\times]t_n,u_n]\in\bar{\mathcal{P}}$,
implicitly defines a measure on $\bar{\mathcal{P}}$ which admits an extension 
$\bar x\in\M(\xi)$ such that $\bar x(F\times\R_+)=(\bar\xi_0-\bar\xi_\infty)(F)
\leq\lambda(F)$ so that $\restr{\bar x}{\bar\F}\in ca(\bar\F)$.
\end{proof}

Corollary \ref{corollary ADM} establishes a general version of the 
Doob Meyer decomposition. In addition it characterises exactly
those processes $\xi$ admitting a countably additive version of such 
decomposition. This characterisation implies a condition hinging on the
uniform countable additivity of the process $\xi$ or, equivalently, a
weak form of countable additivity of the Dol\'{e}ans-Dade measure, namely 
$\restr{\bar x}{\bar\F }\in ca(\bar\F)$. We shall return on this issue
in the following sections.

The existence of Dol\'{e}ans-Dade measures easily translates into that 
of extensions of finitely additive supermartingales, a result 
which may prove useful in problems involving changes of the underlying 
filtration. For $H\subset\OR$ let $H_\omega$ denote the section 
$\{t\in\R_+:(\omega,t)\in H\}$.

\begin{corollary}
\label{corollary supermartingale extension}
Consider an increasing family $(\A_\tau:\tau\in\mathbf{T})$ of algebras of subsets 
of $\Omega$ where $\mathbf{T}\subset 2^{\OR}$ is ordered by 
reverse inclusion
and let $\xi\in\Su$ and $\bar{x}\in\M(\xi)$. There exists a finitely additive 
supermartingale $\xi^*$ on $(\A_\tau:\tau\in\mathbf{T})$ such that $\bar{x}\in\M(\xi^*)$. 
As a consequence 
\begin{enumerate}
\item[(\textit{i})] If $\tau(t)\equiv\Omega\times]t,\infty[\in\mathbf{T}$ and $F\in\A_t\cap\A_{\tau(t)}$ then 
$\xi_{\tau(t)}^*(F) =\xi_t(F)$
\item[(\textit{ii})] If $\tau,\upsilon\in\mathbf{T}$, $F\subset\Omega$ and 
$F_{\tau,\upsilon}\equiv\{\omega\in F:\upsilon_\omega\subset \tau_\omega\}\in\A_\tau\cap\A_\upsilon$ 
then $\xi_\tau^*(F_{\tau,\upsilon})\geq \xi_\upsilon^*(F_{\tau,\upsilon})$
\end{enumerate}
\end{corollary}

\begin{proof}
Fix $\bar{x}\in\M(\xi)$ and define $\xi_\tau^*\in ba(\A_\tau)$ implicitly by letting 
\begin{equation}
\label{Su ext} 
\xi _\tau^*(F)=\bar{\xi}_\infty(F)+\bar{x}((F\times\RR_+)\cap \tau)\qquad F\in\A_\tau
\end{equation}
where $\bar{\xi}_{\infty}$ is an extension of $\xi_\infty$ to 
$2^\Omega$: (\textit{i}) is immediate from (\ref{x}). Given that $\tau\leq \upsilon$ is
equivalent to $\upsilon\subset \tau$ then $F\in\A_\tau$ and $\tau\leq \upsilon$ imply 
$\xi_\tau^*(F)\geq\xi_\upsilon^*(F)$ so that $\xi^*$ is a finitely additive supermartingale 
on $(\A_\tau:\tau\in\mathbf{T})$. Moreover, if $F\subset\Omega$ and 
$F_{\tau,\upsilon}\in\A_\tau\cap\A_\upsilon$ then
$\xi_\tau^*(F_{\tau,\upsilon})\geq \xi_\upsilon^*(F_{\tau,\upsilon})$ is equivalent to
$\bar{x}((F_{\tau,\upsilon}\times\RR_+)\cap \tau)\geq\bar{x}((F_{\tau,\upsilon}\times\RR_+)\cap \upsilon)$ 
which follows from $\bar x$ being positive.
\end{proof}

Whenever $\tau(t)\in\mathbf{T}$ and $\A_t\subset\A_{\tau(t)}$ for all $t\in\R_+$,
Corollary \ref{corollary supermartingale extension} suggests that 
any $\xi\in\Su$ may be consistently extended to any filtration indexed 
by $\mathbf{T}$. Corollary \ref{corollary supermartingale extension}
is an illustration of the importance of finite versus countable additivity.

\section{Two Decompositions}
\label{sec representable}

We shal prove in this section that all finitely additive supermartingales have a 
component that may be represented as a classical supermartingale with respect 
to some $P\in\Prob(\F)$. It should be highlighted that the probability measure 
$P$ involved here emerges endogenously, rather than being given from the outset, 
as in the classical theory. We start with a preliminary result.

\begin{lemma}
\label{lemma measure decomposition}Let $\G\subset\mathcal{H}$ be
two algebras of subsets of $\Omega$ and denote by $ca(\G,\mathcal{H})$ 
and $pfa(\G,\mathcal{H})$ the subspaces of $ba(\G)$ consisting of set functions
which admit a countably additive extension to $\mathcal{H}$ and whose norm
preserving extensions to $\mathcal{H}$ are all purely finitely additive,
respectively. For each $\lambda\in ba(\G)$ there exists a unique way of writing 
\begin{equation}
\lambda =\lambda^e+\lambda^p  \label{measure decomposition}
\end{equation}%
where $\lambda^e\in ca(\G,\mathcal{H})$, $\lambda^p\in pfa(\G,\mathcal{H})$ and 
 $\lambda^e,\lambda^p\geq 0$ if and only if $\lambda\geq 0$. 
\end{lemma}

\begin{proof}
With the aid of the Radon Nikodym theorem it is easily seen that $\lambda
\in ca(\G,\mathcal{H})$ if and only if $\lambda\ll\restr{\bar\lambda}{\G}$ for some 
$\bar\lambda\in ca(\sigma\mathcal{H})$ and, thus, that $ca(\G,\mathcal{H})$ is an 
ideal. Let $\neta{\lambda}$ be an increasing net in $ca(\G,\mathcal{H})_+$ 
with $\bigvee_{\alpha\in\mathbf{A}}\lambda_\alpha=\lambda\in ba(\G)$. Fix 
$\alpha_1\in\mathbf{A}$ arbitrarily and, for given $\alpha_{n-1}$, let 
$\alpha_n\geq\alpha_{n-1}$ be such that $\lambda_{\alpha_n}(\Omega)\geq\lambda(\Omega)-2^{-n}$. 
If $F\in\G$, 
$$
\lambda(F)\geq\lim_n\lambda _{\alpha_n}(F)=
\lambda(\Omega)-\lim_n\lambda _{\alpha_n}(F^c)\geq\lambda(\Omega)-\pi(F^c)=\lambda(F)
$$
But then, $\lambda\ll\sum_n2^{-n}\lambda _{\alpha_n}\in ca(\G,\mathcal{H})$ i.e. 
$\lambda\in ca(\G,\mathcal{H})$. We obtain from Riesz theorem the decomposition 
$ba(\G)=ca(\G,\mathcal{H})+ca(\G,\mathcal{H})^\bot$. The inclusion 
$pfa(\G,\mathcal{H})\subset ca(\G,\mathcal{H})^\bot$ is clear. To prove the converse, 
let $\bar\lambda\in ba(\mathcal{H})$ extend $\lambda\in ca(\G,\mathcal{H})^\bot$. 
Then there exists $G_n\in\G$ such 
that $\abs{\lambda}(G_n^c) +\abs{\bar{\lambda}^c}(G_n)<2^{-n}$. If $G\in\G$
$$
\abs{\lambda(G)}=\lim_n\abs{\lambda(G\cap G_n)}
            =\lim_n\abs{\bar{\lambda}^\bot(G\cap G_n)}
             \leq\abs{\bar{\lambda}^\bot}(G)
$$
i.e. 
$
\norm{\bar{\lambda}^c} +\norm{\bar{\lambda}^\bot}=\norm{\bar{\lambda}}
                                             =\norm{\lambda}\leq\norm{\bar{\lambda}^\bot}
$. 
In other words, $\lambda\in pfa(\G,\mathcal{H})$.
\end{proof}

Lemma \ref{lemma measure decomposition} is a slight generalization of the
classical decomposition of Yosida and Hewitt (by uniqueness the two
decompositions coincide for $\G=\mathcal{H}$). It has though an
important implication here as it implicitly suggests that finitely additive
supermartingales may admit a component that can be represented as a
classical supermartingale with respect to some $P\in\Prob(\F)$.

\begin{proposition}
\label{proposition Su decomposition}
Let $\xi\in\Su_+$. For each $t\in\R_+$ let $\xi_t=\xi_t^e+\xi_t^p$ 
with $\xi_t^e\in ca(\A_t,\F)$ and $\xi_t^p\in pfa(\A_t,\F)$ 
and set $\xi^e=(\xi_t^e:t\in\R_+)$ and $\xi^p=(\xi_t^p:t\in\R_+)$. 
Then%
\begin{equation}
\xi=\xi^e+\xi^p
\label{supermartingale decomposition}
\end{equation}%
is the unique decomposition of $\xi$ such that $\xi^e\in\Su_+$ 
may be represented as a classical $P$ supermartingale $X$ for some 
$P\in\Prob(\F)$ \footnote{The property defined here was called the 
Kolmogoroff property by Bochner \cite[p. 164]{bochner}} while 
$\xi^p$ is positive and orthogonal to all finitely additive 
processes admitting such representation. We say that $\xi^e$ is 
representable and that the pair $(P,X)$ is a representation of $\xi^e$.
\end{proposition}

\begin{proof}
The inclusion $\xi^e\in\Su_+$ was shown in the Introduction. As 
$\xi^p$ is clearly orthogonal to any classical stochastic process,
we only need to prove that $\xi^e$ admits a representation. 
Define the function $T(t) =\norm{\xi _t^e}$ and the set 
$\mathbb{J}=\left\{t\in\R_+ :T(t)>\sup_{u>t}T(u)\right\}$
(with $\sup\varnothing =-\infty $). As $T$ is monotone, $\mathbb{J}$ is
countable; let $\mathbb{C}$ be a countable subset of $\R_+$ such that 
$T[\mathbb{C}]$ is dense in $T[\R_+] $. For each $t\in\R_+$ either 
$t\in\mathbb{J}$ or there is a decreasing sequence $\seq{t}{k}$ in 
$\mathbb{C}$ such that $\lim_kT(t_k)=T(t)$. Let $\seqn{t}$ be an 
explicit enumeration of $\mathbb{D}=\mathbb{C}\cup\mathbb{J}$, choose 
$\bar{\xi}_{t_n}^e\in ca(\F)$ such that $\lrestr{\bar{\xi}_{t_n}^e}{\A_{t_n}}%
=\xi_{t_n}^e$, fix $Q\in\Prob(\F)$ and let $\bar{P}=Q+\sum_n2^{-n}\bar{\xi}_{t_n}^e$
and $P=\norm{\bar{P}}^{-1}\bar{P}$.
Clearly, $P\in\Prob(\F) $, $P\gg\bar{\xi}_{t_n}^e$ for each $n\in\N$. By
construction, for each $t\in\R_+$ and $k>0$ there is 
$t_k\in\mathbb{D}$ such that $t\leq t_k$ 
and $(\xi_t^e-\xi_{t_k}^e)(\Omega)\leq 2^{-k}$. Remark that 
$\lrestr{(\bar{\xi}_t^e-\bar{\xi}_{t_k}^e)}{\sigma\A_t}\in ca(\sigma\A_t)$ 
is the (unique) countably additive extension of 
$\xi_t^{e}-\restr{\xi_{t_{k}}^e} {\A_t}$ to 
$\sigma\A_t$ and is therefore positive. We conclude that 
$\bar{\xi}_t^e(F)=\lim_k\bar{\xi}_{t_{k}}^e(F)$ for each 
$F\in\sigma\A_t$. By Vitali Hahn Saks theorem and its corollaries \cite[III.7.2-3]{bible}, 
$\restr{\bar{\xi}_t^e}{\sigma\A_t}\ll\restr{P}{\sigma\A_t}$,
i.e. $\xi^{e}$ is representable.
\end{proof}

Uniformly countably additive supermartingales play a special role in the 
following section. 
\begin{proposition}
\label{proposition uc decomposition}
Each $\xi\in\Su$ admits a unique decomposition $\xi=\xi^{uc}+\xi^{up}$
where $\xi^{uc}\in\Su^{uc}$ and $\xi^{up}\in\Su^{up}$. 
\end{proposition}

\begin{proof}
Let $\bar{x}\in\M(\xi)$ and let $\bar{x}_{\bar\F}^c$ and $\bar{x}_{\bar\F}^\bot$ be the 
countably and purely finitely additive parts of $\restr{\bar x}{\bar\F}$, respectively. Define 
$\bar{x}'\in\M$ by letting
$$
\bar{x}'(H)=\bar{x}_{\bar\F}^c\left(\Idcond{\bar{x}}{H}{\bar{x}_{\bar\F}^\perp}\right)
\qquad H\subset \OR
$$
Then, by (\ref{take out}), $\bar{x}'_{\bar\F}=\bar{x}_{\bar\F}^c$ 
-- so that $\bar{x}'\in\M^{uc}$. Letting $I_n\in\I{\bar{x}_\F^\perp}$ 
be such that $\bar{x}_\F^c(I_n^c)<2^{-n}$ 
$$
\bar{x}'(H)=\lim_n\bar{x}_{\bar\F}^c\left(I_n\Idcond{\bar{x}}{H}{\bar{x}_{\bar\F}^\bot}\right) 
=\lim_n\bar{x}(I_nH)\leq\bar{x}(H)\qquad H\subset \OR
$$
Clearly, $\bar x''=\bar x-\bar x'\in\M^{up}$. Thus the set
$\M^*(\xi)=\{\bar{y}\in\M^{uc}:\bar{y}\leq\bar{x}\text{ for some }\bar{x}\in\M(\xi)\}$
is non empty and, we claim, it admits a maximal element with respect to the
partial order $\geq_{\bar\F}$ defined by letting $\bar y\geq_{\bar\F}\bar y'$ whenever 
$\bar{y}_{\bar\F}\geq\bar{y}'_{\bar\F}$. In order to apply Zorn 
lemma, consider an increasing net $\nneta{\bar y^\alpha}$ in $\M^*(\xi)$ and let 
$\bar x^\alpha\in\M(\xi)$ be such that $\bar y^\alpha\leq\bar x^\alpha$ 
for all $\alpha\in\mathbf{A}$. Define 
$\bar x^\circ(H)=\text{LIM}_{\alpha\in\mathbf{A}}\bar x^\alpha(H)$, 
$H\subset \OR$ -- where $\text{LIM}$
denotes the Banach limit functional introduced in \cite{agnew morse}. By
linearity, $\bar x^\circ\in\M(\xi)$. The inequality 
$\bar x^\circ(H)\geq\liminf_{\alpha\in\mathbf{A}}\bar x^\alpha(H)
\geq \liminf_{\alpha\in\mathbf{A}}\bar y^\alpha(H)$ which holds for any 
$H\subset\OR$ implies that 
$\bar x^\circ\geq_{\bar\F}\bar y^\alpha$ for all $\alpha\in\mathbf{A}$
i.e. that $\bar x^\circ$ is an upper bound for $\nneta{\bar y^\alpha}$. 
Let $\bar{x}^{uc}$ be a maximal element of $\M^*(\xi)$, 
let $\bar{x}^*\in\M(\xi)$ be such that $\bar{x}^{uc}\leq\bar{x}^*$ and define 
$\bar{x}^{up}=\bar{x}^*-\bar{x}^{uc}\in\M$. Let $\xi^{uc},\xi^{up}\in\Su$ be uniquely 
defined by the condition $\bar{x}^{uc}\in\M(\xi^{uc})$, $\bar{x}^{up}\in\M(\xi^{up})$, 
$\xi_\infty^{uc}=\xi_\infty^c$ and $\xi_\infty^{up}=\xi_\infty^\bot$. 
By construction, $\xi^{uc}\in\Su^{uc}$. Decompose $\bar{y}\in\M(\xi^{up})$ as 
$\bar{y}'+\bar{y}''$ where $\bar{y}'\in\M^{uc}$
and $\bar{y}''\in\M^{up}$, as in the first step of
this proof. From $\bar{x}^{uc}+\bar{y}'\leq\bar{x}^{uc}+\bar{y}\in\M(\xi)$ 
and the fact that $\bar{x}^{uc}$ is $\geq_{\bar\F}$ maximal, we deduce $\bar{y}'=0$ or, 
equivalently, $\xi^{up}\in\Su^{up}$. If $\xi=\kappa^{uc}+\kappa^{up}$ 
were another such decomposition, and $k^{up}$ and $k^{uc}$ the associated 
Dol\'{e}ans-Dade measures, then from $k^{up}\leq x^{uc}+x^{up}$ and Hahn Banach
one may find $\bar{k}^{up}\in\M(\kappa^{up})$ such that 
$\bar{k}^{up}\leq\bar{x}^{uc}+\bar{x}^{up}$. However, since 
$\bar{k}^{up}\perp\bar{x}^{uc}$, this implies $\bar{k}^{up}\leq\bar{x}^{up}$ 
while the converse is obtained mirrorwise. In other words $\kappa^{up}$ 
and $\xi^{up}$ induce the same Dol\'{e}ans-Dade measure; in addition,
$\kappa_\infty^{up}=\xi_\infty^{up}=\xi_\infty^\perp$. The claim follows
from Theorem \ref{theorem extension}(\textit{iii}).
\end{proof}

\section{Increasing Processes}
\label{sec integration}

Fix $P\in\Prob(\F)$ and let $\AP$ denote the set of processes 
$(A_t:t\in\R_+)$ such that $A_\infty\in L^1(P)$ and $P(0=A_0\leq A_t\leq A_u)=1$
for each $0\leq t\leq u<\infty$. Of course, if $A\in\AP$ and $A'$ is a 
modification of $A$ (i.e. $P(A'_t=A_t)=1$ for all $t\in\R_+$) then $A'\in\AP$.
Put $\mathfrak{A}=\bigcup_{P\in\Prob(\F)}\AP$.

\begin{lemma}
\label{lemma increasing}
Let $A\in\AP$. Then there is $F\in\F$ with $P(F^c) =0$ and a modification $A'$ of 
$A$ such that for each $t\leq u$, $0=A'_0\leq A_t'\leq A_u'$ on $F$. If in addition 
$P(A_t)=\lim_nP(A_{t+2^{-n}})$ then $A'$ and $F$ may be chosen to be right 
continuous at each $t\in\R_+$ and for each $\omega\in F$.
\end{lemma}

\begin{proof}
As in the proof of Proposition \ref{proposition Su decomposition}, there exists
a countable subset $\mathbb{D}$ of $\R_+$ with the property that for each $t\in\R_+$ 
and $\epsilon>0$ there is $d\in\mathbb{D}$ such that $d\geq t$ and $P(A_t)>P(A_d)-\epsilon$. 
Define $F=\bigcap_{\{d,d'\in\mathbb{D}:d>d'\} }\{A_d\geq A_{d'}\}$: clearly, $P(F^c)=0$.
Let $\mathbb{D}(t)=\{d\in\mathbb{D}:d\geq t\}$ and 
$A_t'=\inf_{d\in\mathbb{D}(t)}A_d$. By definition of 
$\mathbb{D}$, $A'_t\geq A_t$ but $P(A_t')=P(A_t)$ so that $P(A_t=A'_t)=1$.
If $A$ is right continuous in mean the same conclusion holds even if we replace
$\mathbb{D}(t)$ with $\mathbb{D}^+(t)=\{d\in\mathbb{D}:d>t\}$. However
$A'$ is right continuous on $F$ since $\mathbb{D}^+(t)=\bigcap_{u>t} \mathbb{D}^+(u)$.
\end{proof}

For $H=(F_0\times\{0\})\cup\bigcup_{n=1}^N(F_n\times]t_n,u_n])%
\in\mathcal{\bar{P}}$ and $A\in\AP$ the integral 
$\int \set{H}dA$ has an obvious definition, namely 
$\sum_{n=1}^N\set{F_n}(A_{u_n}-A_{t_n})$. In the following Theorem we obtain 
an extension of this integral together with a characterization of increasing 
processes in terms of their Dol\'{e}ans-Dade measure.

\begin{theorem}
\label{theorem stochastic integral}
Let $\bar{x}\in\M$. The following are equivalent:
\begin{enumerate}
\item[(\textit{i}). ] There exists $P\in\Prob(\F)$ and, given $P$, a unique
(up to modification) $A\in\AP$ such that 
\begin{equation}
\bar{x}(H)=P\int \set{H}dA\qquad H\in\mathcal{\bar{P}}
\label{iso}
\end{equation}
\item[(\textit{ii)}. ] $\bar{x}\in\M^{uc}$ and $\bar{x}(\{0\})=0$;
\item[(\textit{iii}). ] There exists $P\in\Prob(\F)$ such that for each 
$h\in L^1(\bar{x})$ the equation%
\begin{equation}
\bar{x}(bh)=P(bI_{\bar{x}}(h))\qquad b\in\B(\F)   
\label{isomorphism}
\end{equation}%
admits a unique solution $I_{\bar{x}}(h)\in L^1(P)$ such that $I_{\bar{x}}(\set{\{0\}})=0$.
\end{enumerate}
$I_{\bar{x}}:L^1(\bar{x})\rightarrow L^1(P)$ as defined in (\ref{isomorphism})
is a positive, continuous, linear functional such that $\norm{I_{\bar x}} =\norm{\bar x}$
and that $\lim_nI_{\bar{x}}(h_n)=0$ whenever $\sup_n\abs{h_n}\in L^1(\bar{x})$ 
and $\lim_nP^*(h^*_n>\eta)=0$ for each $\eta >0$ -- where 
$h^*_n\equiv\sup_t\abs{h_{n,t}}$ and $P^*$ is the outer measure generated by $P$. 
\end{theorem}

\begin{proof}
Let us start remarking that one may easily identify $\bar\F$ with $\F$,
as we shall now do. Under (\textit{i}), $\bar{x}(\{0\})=P\int\set{\{0\}} dA=0$ 
and $\bar{x}(F) =\bar{x}(F\times]0,\infty[)=P(\set{F}A_\infty)$ 
for any $F\in\F$. Assume (\textit{ii}) and fix 
$P=(\norm{Q} +\norm{\bar{x}})^{-1}(Q+\restr{\bar x}{\bar\F})$ 
for some $Q\in\Prob(\F)$. By Theorem \ref{theorem conditioning}, 
for each $h\in L^1(\bar x)$ we may define
$$
I_{\bar{x}}(h)=\cond{\bar{x}}{h}{\bar\F}\frac{d\restr{\bar x}{\bar\F}}{dP}\in L^1(P) 
$$
By (\ref{conditioning}), $I_{\bar{x}}(h)$ is a solution to (\ref{isomorphism}); 
moreover the operator $I_{\bar{x}}$ is positive, linear and has norm $\norm{\bar{x}}$, 
by Theorem \ref{theorem conditioning}; $\bar{x}(\{0\}) =0$ implies 
$I_{\bar{x}}(\set{\{0\}})=0$, $P$ a.s.. Any other solution $J(h)\in L^1(P)$ to 
(\ref{isomorphism}) satisfies $P(bJ(h))=P(bI_{\bar{x}}(h))$ for all 
$b\in L^\infty(P)$ i.e. $P(J(h)=I_{\bar{x}}(h))=1$. 
Assume (\textit{iii}) and define $A_t=I_{\bar{x}}(\set{]0,t]})$, 
$A=(A_t:t\in\R_+)$ and let 
$H=(F_0\times\{0\})\cup\bigcup_{n=1}^N(F_n\times]t_n,u_n])\in\mathcal{\bar{P}}$. 
Then $A\in\AP$ and, up to a $P$ null set
\begin{eqnarray*}
I_{\bar{x}}(\set{H})%
&=&\sum_{n=1}^N\rcond{\bar x}{\set{F_n\times]t_n,u_n]}}{\bar\F}\frac{d\restr{\bar x}{\bar\F}}{dP}  \\
&=&\sum_{n=1}^N\set{F_n}\rcond{\bar x} {\set{]t_n,u_n]}} {\bar\F}\frac{d\restr{\bar x}{\bar\F}}{dP}  \\
&=&\sum_{n=1}^N\set{F_n}(A_{u_n}-A_{t_n})  \\
&=&\int \set{H}dA 
\end{eqnarray*}
But then (\ref{iso}) follows from (\ref{isomorphism}). If $B\in\AP$ 
also meets (\ref{iso}) then for $h=\set{F\times]0,t]}$ and $F\in\F$ 
we conclude that $P(FA_t)=P(FB_t)$ from which we deduce uniqueness. 
It is clear from (\ref{iso}) that $I_{\bar x}$ is linear, positive 
and that $\norm{I_{\bar x}} =\norm{\bar x}$.

If $\seqn{h}$ is a sequence in $L^1(\bar{x})$ with the above properties 
then so is $\sseqn{\abs{h_n}}$. Given that $I_{\bar x}$ is positive, it 
is enough to prove the claim for $h_n\geq 0$. Observe that 
$\bar{x}(h_n\geq\eta)\leq\bar{x}(h^*_n\geq\eta)$; moreover, 
$\restr{\bar x}{\bar\F}\ll P$ implies that, in restriction to 
$2^\Omega\otimes\{\varnothing,\R_+\}$, $\bar x\ll P^*$. But then, $h_n$ 
converges to $0$ in $\bar{x}$ measure and, by \cite[theorem III.3.6]{bible}, 
in $L^1(\bar{x})$. Given (\ref{isomorphism}) this is equivalent to 
$I_{\bar{x}}(h_n)$ converging to $0$ in $L^1(P)$.
\end{proof}

The equivalence of (\textit{i}) and (\textit{ii}) establishes a correspondence 
between $\mathfrak{A}$ and $\M^{uc}$ which compares to the classical 
(and well known) characterization of increasing processes as measures 
given by Meyer \cite[VI.65, p. 128] {dellacherie meyer} (see also \cite%
[p. 6]{meyer cours}). Meyer's result, which ultimately delivers the Doob 
Meyer decomposition, focuses however on countable additivity over 
$\F\otimes\mathcal{B}(\R_+)$; we rather require this property relatively 
to $\bar\F$. In Theorem \ref{theorem Doob Meyer} below we show that indeed 
this is enough to obtain a suitable version of Doob Meyer decomposition. 
Although there are connections between these two properties, it is noticeable 
that the latter is independent of the given filtration. On should also remark
that we do not assume the existence of an underlying probability $P\in\Prob(\F)$
but rather deduce it.

Each $\bar x\in\M$ may be considered in restriction to special classes
of functions such as the set $C$ of functions $f:\OR\rightarrow\R$
with continuous sample paths and bounded support (i.e. such that $f(t)=0$ 
for all $t$ larger than some $T$). Let $\mathcal{C}$ be the $\sigma$ algebra 
on $\OR$ generated by $C$. 
\begin{lemma}
\label{lemma Dini Daniell}
Let $\bar x\in\M^c$. There exist $\alpha^c\in ca(\mathcal{C})_+$, $P\in\Prob(\F)$ 
and $A^c\in\AP$ right continuous such that
\begin{equation}
\label{compensator}
\bar x(f)=\alpha^c(f)=P\int fdA^c\qquad f\in L(\bar x)\cap C
\end{equation}
\end{lemma}
\begin{proof}
In order to apply Daniell theorem, consider a sequence $\seqn{h}$ 
in the vector lattice $L(\bar x)\cap C$ decreasing to $0$ and fix 
$T$ such that $\bar x(\abs{h_1-h^T_1})<\epsilon$, where $h_n^T=h_n\set{]0,T]}$. 
Let $h_n^{T,*}=\sup_th^T_n(t)$. A simple application of Dini's theorem 
for each $\omega\in\Omega$ guarantees that the sequence $\seqn{h^{T,*}}$ 
converges to $0$ pointwise; moreover, by continuity of the sample 
paths, $h^{T,*}_n$ is in fact $\F$ measurable. Thus Theorem 
\ref{theorem stochastic integral} implies that 
$
\lim_n\bar x(h_n)\leq\epsilon+\lim_n\bar x(h^T_n)
=\epsilon+\lim_nP(I_{\bar x}(h^T_n))=\epsilon
$.
In other words, the restriction of $\bar x$ to $C$ is a Daniell integral 
and as such it admits the representation as the integral with respect to 
some $\alpha^c\in ca(\mathcal{C})$. Observe that $F\times]t,\infty[\in\C$ 
for all $F\in\F$ and $t\in\R_+$. Fix $P\in\Prob(\F)$ as in Theorem 
\ref{theorem stochastic integral}(\textit{i}) and define $\alpha^c_t\in ba(\F)$ 
as $\alpha^c(F\times]0,t])$ for each $F\in\F$. Since 
$\alpha^c_t\leq\bar x_{\bar\F}\ll P$, denote by $A^c_t$ the Radon Nikodym 
derivative of $\alpha^c_t$ with respect to $P$. We deduce that 
$P((A^c_u-A^c_t)\set{F})=\alpha^c(F\times]t,u])=\bar x(F\times]t,u])\geq 0$, 
so that $A^c\in\AP$, and that 
$\lim_nP(\set{F}(A^c_{t+2^{-n}}-A^c_t))=\lim\alpha^c(F\times]t,t+2^{-n}])=0$ 
(by countable additivity) for each $F\in\F$ so that $A^c_t=\lim_nA^c_{t+2^{-n}}$ 
up to a $P$ null set. By Lemma \ref{lemma increasing}, we obtain that $A^c$ 
admits a modification which is right continuous.
\end{proof}

\begin{theorem}
\label{theorem Doob Meyer}
Let $\xi\in\Su$. Then $\xi\in\Su^{uc}$ if and only if there exist $P\in\Prob(\F)$, 
$M\in L^1(P)$ and $A^p\in\AP$ which is adapted, right continuous in mean and such 
that
\begin{equation}
\label{Doob Meyer}
\xi_t(F)=P(\set{F}(M-A^p_t))\qquad t\in\R_+\ \text{and}\ F\in\F_t
\end{equation}
and that 
\begin{equation}
\label{natural}
P\left(b\int hdA^p\right)=P\int M(b)_-hdA^p\qquad b\in L^\infty(P),\ h\in\B(\sigma\mathcal{P})
\end{equation}
where $M(b)=(\condP{b}{\F_t}:t\in\R_+)$.
\end{theorem}
\begin{proof}
We use the notation of Lemma \ref{lemma Dini Daniell} and the inclusion
$\sigma\mathcal{P}\subset\C$. Let $d=\{t_1\leq t_2\leq\ldots\leq t_N\}$ be 
a finite sequence in $\R_+$ and  define
\begin{equation}
\label{Pd projection}
\Pd(f)=\sum_{n=1}^{N-1}\cond{P}{f_{t_n}}{\F_{t_n}}\set{]t_n,t_{n+1}]}  \quad
\text{and} \quad \bar x^d(f)=\bar x(\Pd(f))\qquad f:\R_+\rightarrow L^1(P)
\end{equation}
Denote by $\alpha^d$ the restriction of $\bar x^d$ to $\F\otimes 2^{\R_+}$. 
On the one hand it is easily seen that $\bar x^d_{\bar\F}\ll P$ so that, as 
in the proof of Theorem \ref{theorem stochastic integral}, we can associate 
to $\bar x^d$ a process $A^d\in\AP$, by letting 
$A^d_tdP=\cond{\bar x^d}{\set{]0,t]}}{\bar\F}d\bar x^d_{\bar\F}$.  
On the other hand, (\ref{compensator}) implies
$
P\int fdA^d=\bar x(\Pd(f))=\alpha^c(\Pd(f))
$.
Consider the case in which $f=bh$ where $b\in L^\infty(P)$ and 
$h$ is bounded, adapted and left continuous. 
Let $d_n=\{k2^{-n}:k=0,\ldots,2^{2n}\}$ and observe that, by \cite[VI.2, p. 67]
{dellacherie meyer}, there exists a $P$ null set $F\in\F$ outside of which 
$\lim_n\Pn(b)_t=M(b)_{t-}$ and $\lim_n\Pn(h)_t=h_t$ for each $t\in\R_+$
(as $h$ is left continuous and adapted). Given that $\alpha^c$ is countably 
additive in restriction to $\mathcal{C}$, we conclude
\begin{equation}
\label{comp} 
\lim_nP\left(b\int hdA^{d_n}\right)=\lim_n\alpha^c(\Pn(bh))
=\lim_n\alpha^c(\Pn(b)\Pn(h))=\alpha^c(M(b)_-h)
\end{equation}
Define then $\alpha^p\in ba(\F\otimes 2^{\R_+})$ implicitly
as $\alpha^p(H)=\alpha^c(M(\set{H})_-)$. 
Then from (\ref{comp}) we deduce that $\sseqn{\alpha^{d_n}}$ converges to $\alpha^p$
and, by \cite[III.7.3, p. 159]{bible}, that $\alpha^p_{\bar\F}\ll P$. Let $A^p\in\AP$ 
be the increasing process associated to $\alpha^p$. Thus for
every bounded, adapted and left continuous process $h$ and every $b\in L^\infty(P)$
we have
$$
P\left(b\int hdA^p\right)=\bar x(M(b)_-h)=\alpha^c(M(b)_-h)=P\int M(b)_-hdA^p
$$
which delivers (\ref{Doob Meyer}) if we only let $M=d\xi_\infty/dP+A^p_\infty$,
$b=\set{F}$ with $F\in\F_t$ and $h=\set{]t,\infty[}$. In addition, if 
$F\in\F$, $s\leq t$  and $h_{F,t}=\set{F}-\condP{\set{F}}{\F_t}$, then $M(h_{F,t})_{s-}=0$
so that 
$$
P(h_{F,t}A^p_t)=\alpha^p(h_{F,t}\set{[0,t]})=P\int_0^tM(h_{F,t})_-dA^p=0
$$
Therefore, replacing $A^p_t$ with $\condP{A^p_t}{\F_t}$, we may assume that 
$A^p\in\AP$ is adapted. Eventually, letting $h_n=\set{]t,t+2^{-n}]}$ 
we conclude that $0=\lim_n\alpha^c(h_n)=\lim_nP(A^p_{t+2^{-n}}-A^p_t)$
and, thus, that $A^p$ is right continuous in mean. 
That (\ref{Doob Meyer}) implies $\xi\in\Su^{uc}$ is obvious.
\end{proof}

With a complete filtration Theorem \ref{theorem Doob Meyer} implies
that $A^p$ may be chosen to be adapted and right continuous.

We want to emphasize that the existence of the decomposition (\ref{Doob Meyer}) does 
not depend on the underlying filtration. 

\begin{corollary}
\label{corollary enlargement}
Let $\xi\in\Su$ and let $\mathbb{D}\subset\R_+$ be such that 
$\xi_t=\sup_{d\in\mathbb{D}(t)}\restr{\xi_d}{\A_t}$
where $\mathbb{D}(t)=\{d\in\mathbb{D}:d\geq t\}$. Then, $\xi$ admits a Doob 
Meyer decomposition if and only if $\xi^{\mathbb{D}}=(\xi_d:d\in\mathbb{D})$ does.
\end{corollary}

\begin{proof}
Given that, by Theorem \ref{theorem Doob Meyer}, the Doob Meyer decomposition 
is equivalent to $\xi\in\Su^{uc}$, the direct implication is obvious. 
As for the converse, choose $\bar x^{\mathbb{D}}\in\M(\xi^{\mathbb{D}})$
to be countably additive in restriction to $\F\otimes\{\varnothing,\mathbb{D}\}$. 
If $t\in\R_+$ and $F\in\A_t$ then 
\begin{equation}
\label{dominus} 
\abs{\xi_t}(F)\leq\sup_{d\in\mathbb{D}(t)}\abs{\xi_d}(F)
\leq\abs{\xi_\infty}(F)+\bar x^{\mathbb{D}}(F\times\mathbb{D})\equiv\lambda(F)
\end{equation}
The claim then follows from Corollary \ref{corollary ADM}. 
\end{proof}

Corollary \ref{corollary enlargement} makes clear that decomposition (\ref{Doob Meyer}) 
is a property that involes any subset $\mathbb{D}$ which is dense for the range of $\xi$
and we know from the proof of Proposition \ref{proposition Su decomposition} that 
this may be taken to be countable. The class $D$ property may thus be replaced by a
corresponding property, the class $D_\sigma$, in which the stopping times are restricted
to have countable range, see \cite{io STAPRO}.

\section{The Bichteler Dellacherie Theorem without Probability}
\label{sec bichteler}

Let $f:\OR\rightarrow\R$ be adapted to the filtration, define 
$f^*=\sup_{t\in\R_+}\abs{f_t}$ and let $\F$ be such that $f^*$ is $\F$ measurable.
The starting point of this section are the sets%
\begin{equation}
\K=\left\{\int hdf:h\ \text{is}\ \mathcal{P}\ \text{simple} ,\abs{h}\leq 1\right\}
\quad\text{and}\quad\C=\K -\B(\F)_+  
\label{K}
\end{equation}%
Bichteler and Dellacherie start from the assumption that $\K$ is bounded in $L^0(P)$ 
for some given $P\in\Prob(\F)$ and that $f$ is right continuous with left limits outside 
some $P$ null set. These two properties are then shown to imply that for 
given $\eta>0$ there exists $\delta>0$ such that $d\set{F}\notin\C\set{F}$ for all $F\in\F$ 
such that $P(F)>\eta$. We take inspiration from this separating condition to define a 
concept of boundedness suitable for our setting. To this end we denote by $\U$ a collection 
of subets of $\Omega$ with the following properties:
\begin{assumption}
\label{Ass} 
There exists $\lambda_0>0$ such that 
\begin{equation}
\{\lambda\set{U}:\lambda\geq\lambda_0\}\cap\C\set{U}=\varnothing\qquad U\in\U
\label{separation}
\end{equation}
Moreover, $U,V\in\U$ imply $U\subset\{f^*<n\}$ for some $n$ and $U\cup V\in\U$. 
\end{assumption}
A violation of (\ref{separation}) indicates that the set $\K$ is unbounded 
relatively to some $U\in\U$. Both sets in (\ref{separation}) are convex subsets of $\B(2^U)$ 
and $\set{U}\C$ contains $-\set{U}$ as an internal point. By ordinary properties of the support 
functional \cite[lemma V.I.8(f), p. 411]{bible}, the Hahn Banach theorem and 
\cite[lemma V.II.7, p. 417]{bible} we conclude that for each $U\in\U$ there is 
$\hat{m}_U\in ba(2^U)$ such that $\sup_{x\in\set{U}\C}\hat{m}_U(x)\leq 1=\lambda_0\hat{m}_U(U)$. 
The inclusion $-\B(2^U)_+\subset \set{U}\C$ implies that $\hat{m}_U\geq 0$. 
By defining $m_U\in ba_+$ implicitly by 
$$
m_U(F)=\dfrac{\hat{m}_U(F\cap U)}{\hat{m}_U(U)} \qquad F\subset\Omega
$$
we have completed the proof of the following:
\begin{lemma}
\label{lemma yan}
Let $f:\OR\rightarrow\R_+$
satisfy Assumption \ref{Ass} and define the set 
\begin{equation}
\mathcal{M} =\left\{ m\in ba_+:\norm{m}=1,\ \sup_{x\in\C }m(x)\leq \lambda_0\right\}
\label{M(C)}
\end{equation}%
For each $U\in\U$ there exists $m_U\in\mathcal{M}$ such that $m_U(U)=1$.
\end{lemma}
Fix now $m\in\mathcal{M}$ (so that $f^*\in L^1(m)$) and let 
$\xi^e$ and $\xi^p$ be the components of the finitely additive supermartingale 
$(\restr{m}{\F_t}:t\in\R_+)$ as of (\ref{supermartingale decomposition}). 
Set also $\I{t}=\I{\xi_t^p}$ (see (\ref{I})), let $(P,X)$ be a 
representation for $\xi^e$ and observe that $-\xi^p\in\Su$ and that 
$\M(\xi^e)=\M(-\xi^p)$. Fix an extension $\bar\xi_\infty\in ba(\F)$ of $\xi_\infty$ 
to $\F$ and $\bar{x}\in\M(\xi^e)$ and define
$$
\bar{\xi}_t^p(F) =\bar\xi_\infty^p(F)-\bar x(F\times]t,\infty[)\qquad F\in\F
$$
The collection $(\bar\xi_t^p:t\in\R_+)$ is then increasing with $t$. 
For each $b\in\B(\F)$ and $u\geq t$, $F\in\I{t}$ implies 
$\xi_\infty^p(b\set{F})=\bar{x}(b\set{F\times]t,\infty[})$ and thus
\begin{equation}
\bar\xi_u^p(b\set{F})=\xi_\infty^p(b\set{F})-\bar{x}(b\set{F\times]u ,\infty[})
=\bar{x}(b\set{F\times]t,u]})
\label{fund}
\end{equation}%
Let now $d=\{t_1\leq\ldots\leq t_N\}$, 
\begin{equation}
\Pd=\left\{F_0\times\{0\}\cup\bigcup_{n=1}^{N-1}F_n\times]t_n,t_{n+1}] :
F_0\in\F_0,F_n\in\F_{t_n},1\leq n\leq N-1\right\}  
\label{Pd}
\end{equation}
and choose $F_n\in\I{t_n}$ $1\leq n<N$ 
and set
$$
F^d=\bigcup_{n=1}^{N-1}F_n\times]t_n,t_{n+1}]\quad \text{and}
\quad f^d=f_0\set{\{0\}}+\sum_{n=1}^{N-1}f_{t_{n+1}}\set{]t_n,t_{n+1}]}
$$ 
By (\ref{fund})%
\begin{eqnarray*}
\sum_{n=1}^{N-1}\xi_{t_{n+1}}^p((f_{t_{n+1}}-f_{t_n})\set{F_n})
&=&\sum_{n=1}^{N-1}\bar{x}((f_{t_{n+1}}-f_{t_n})\set{F_n\times]t_n,t_{n+1}]}) \\
&=&\bar{x}(f^d\set{F^d})
-x\left(\sum_{n=1}^{N-1}f_{t_n}\set{F_n\times]t_n,t_{n+1}]}\right) \\
&=&\bar{x}(f^d\set{F^d})+P\sum_{n=1}^{N-1}f_{t_n}\set{F_n}(X_{t_{n+1}}-X_{t_n})
\end{eqnarray*}%
i.e. 
\begin{eqnarray}
\notag
m\left(\int \set{F^d}df\right)&=&\sum_{n=1}^{N-1}m((f_{t_{n+1}}-f_{t_n})\set{F_n})\\ 
&=&\sum_{n=1}^{N-1}(\xi_{t_{n+1}}^p+\xi_{t_{n+1}}^e)((f_{t_{n+1}}-f_{t_n})\set{F_n})\\ \notag
&=&\bar{x}(f^d\set{F^d})+P\sum_{n=1}^{N-1}\set{F_n}(f_{t_{n+1}}X_{t_{n+1}}-f_{t_n}X_{t_n}) 
\label{parts}
\end{eqnarray}%
Assume that $H=H_0\set{\{0\}} +\bigcup_{n=1}^{N-1}H_n\set{]t_n,t_{n+1}]}\in\Pd$. 
Then by (\ref{x}) 
$$
\bar{x}(H) =P\left\{\set{H_0}(X_{0}-X_{\infty})+\sum_{n=1}^{N-1}\set{H_n}(X_{t_n}-X_{t_{n+1}})\right\}
$$
i.e. $\restr{\bar{x}}{\Pd}$ is countably additive. 
Replacing $F^d$ with a sequence $\sseq{F^{d,k}}{k}$ such that 
$\lim_{k}P\left(\bigcap_{n=1}^{N-1}F_{n-1}^{k}\right) =1$ we thus deduce then from 
(\ref{parts}) 
\begin{equation}
\lim_k m\left(\int \set{F^{d,k}}df\right) =\bar{x}(f^d)+P(f_\infty X_\infty-f_0X_0)  
\label{int}
\end{equation}%
Replace $f$ with $\int \set{F\times]t,u]} df$ where $t\leq u$ and $F\in\F_t$ and choose 
$d$ such that $F\times]t,u]\in\Pd$. We also deduce%
\begin{equation}
\lim_{k}m\left(\set{F}\int_t^u \set{F^{d,k}}df\right)
=P(F(f_uX_u-f_tX_t))+\bar{x}(f^d\set{F\times]t,u]})  
\label{local}
\end{equation}

\begin{theorem}
\label{theorem Bichteler}
Let $f\in\R^{\OR }$ satisfy Assumption \ref{Ass}. Then there exists $P\in\Prob(\F)$ 
and a $P$ positive supermartingale $X$ such that $Xf$ is a $P$ quasimartingale. If 
there is $Q\in\Prob(\F)$ and $\eta>0$ such that $Q(f^*<\infty)=1$ and that 
$F\in\F$ and $Q(F\cap\{f^*<k\})\geq\eta$ imply $F\cap\{f^*<k\}\in\U$ then 
for any $\delta>\eta$ the pair $(P,X)$ above can be chosen such that $P(X_\infty=0)<\delta$.
\end{theorem}

\begin{proof}
By Lemma \ref{lemma yan} for fixed $n>n_0$ there is $m\in\mathcal{M}$ 
such that $m(f^*>n)=0$ so that $m(f^*)<\infty$. By (\ref{fund}), 
$$
\abs{\bar{x}(f^d)}\leq\sum_{n=1}^{N-1}(\bar\xi_{t_{n+1}}^p-\bar\xi_{t_n}^p)
\left(\sup_{1<j\leq N}\abs{f_{t_{j}}}\right) 
=\bar\xi_\infty^p\left(\sup_{1<j\leq N}\abs{f_{t_{j}}}\right)\leq m(f^*)
$$
Let $h_n^d$ be the sign of $\cond{P}{f_{t_{n+1}}X_{t
_{n+1}}}{\F_{t_n}}-f_{t_n}X_{t_n}$ and $h^d=\sum_nh_n^d\set{]t_n,t_{n+1}]}$. 
By (\ref{local})%
\begin{eqnarray*}
P\sum_n\dabs{\cond{P}{f_{t_{n+1}}X_{t_{n+1}}}{\F_{t_{n}}}-f_{t_n}X_{t_n}} 
&=&P\sum_nh_n^d(f_{t_{n+1}}X_{t_{n+1}}-f_{t_n}X_{t_n}) \\
&=&\lim_{k}m\left(\int F^{d,k}h^{d}df\right) -\bar{x}(f^dh^d) \\
&\leq &\sup_{k\in\K}m(k)+m(f^*)
\end{eqnarray*}%
Since such bound is uniform in $d$ this proves the first claim.

Let $Q\in\Prob(\F)$ and $\eta$ be as in the claim. In search of a contradiction,
assume that for some $\delta>\eta$ and all representations $(P,X)$ of $\xi^e$ we 
have $P(X_\infty=0)>\delta$. Fix $\epsilon=\delta-\eta$, define $P^\epsilon=%
(1-\epsilon)Q+\epsilon P$ and let $(P^\epsilon,X^\epsilon)$ be the corresponding
representation. Given that $P^\epsilon$ and $\xi_\infty^p$ are orthogonal 
and that $\F$ is a $\sigma$ algebra, there is an $\F$ measurable subset $F$ of 
$\{X_\infty^\epsilon=0\}$ such that $\xi_\infty^p(F)=m(F)=0$ and 
$P^\epsilon(F)>\delta$ so that 
$Q(F\{f^*<k\})\geq\frac{\delta-\epsilon}{1-\epsilon}=\frac{\eta}{1-\epsilon}$ 
for some integer $k$. We conclude that for each $m\in\mathcal{M}$ there exists 
$h_m\in\B(\F)$ such that $0\leq h_m\leq 1$, $Q(h_m)\geq\frac{\eta}{1-\epsilon}$, 
$\{h_m>0\}\in\U$ and $m(h_m)=0$. Denote by $H$ the corresponding collection. Then
$$
\sup_{m\in\mathcal{M}}\inf_{h\in H}m(h)=0
$$
Endow $\B(\F)$ with the norm topology and $ba(\F)$ with the weak$^*$ topology.
One easily remarks that $H$ and $\mathcal{M}$ are convex sets (as $\U$ is closed 
with respect to unions) and that $\mathcal{M}$, being a closed subset of the unit sphere 
of $ba(\F)$, is compact; moreover, the function $(m,h)\rightarrow m(h)$ 
is bilinear and separately continuous. Sion's \cite[corollary 3.3, p. 174]{sion} 
version of the minimax theorem therefore applies, yielding the conclusion%
$$
\inf_{h\in H}\sup_{m\in\mathcal{M}}m(h)=0
$$
There is then a sequence $\seqn{h}$ in $H$ such that $\sup_{m\in\mathcal{M}}m(h_n) <2^{-n}$. 
Given that $H$ is convex we may equivalently replace $h_n$ by a convex combination 
$\sum_{j=0}^{J}\alpha_jh_{n+j}$. As a consequence of Komlos lemma \cite[theorem 1, p. 218]{komlos} 
(but see also \cite[theorem 6, p. 184]{schwartz}) there is no loss of generality in assuming that 
the sequence $\seqn{h}$ converges $Q$ a.s. to some $h'$. By Egoroff theorem, we can choose $F\in\F$ 
such that $h_n$ converges uniformly to $h'$ on $F$ and that, letting $h=h'\set{F}$,
$$
Q(h)=\lim_nQ(h_n\set{F})\geq\lim_nQ(h_n)-Q(F^c)\geq\frac{\eta}{1-\epsilon/2}>\eta
$$
while%
$$
\sup_{m\in\mathcal{M}}m(h)\leq\lim_n\sup_{m\in\mathcal{M}}m(h_n) =0
$$
Both inequalities remain true if we replace $h$ by $U=\{h>a;f^*<1/a\}$ 
for $a$ sufficiently small. Then $U\in\U$ but $m(U)=0$ for all $m\in\mathcal{M}$ so that 
Lemma \ref{lemma yan} fails, a contradiction.
\end{proof}

\end{document}